\documentclass[10pt]{article}
\usepackage{latexsym, amscd, amsmath, graphics, rotating, array}
\usepackage{amssymb,amsmath,multibox}
\makeindex
\title{Functional equations for local normal \\zeta functions of nilpotent groups}
\author{Christopher Voll\medskip\\  with an appendix by \\ Arnaud Beauville}

\newenvironment{acknowledgements}{\bigskip\noindent\textsl{Acknowledgements.}\rm}

\newenvironment{proof}{\noindent\textbf{Proof.}\rm}{\hfill$\Box$ \\ \hspace{0.1cm}}

\newenvironment{proofnodot}{\bigskip\noindent\textbf{Proof}\rm}{\hfill$\Box$ \\ \hspace{0.1cm}}

\newtheorem{lemma}{Lemma}
\newtheorem{assumption}{Assumption}
\newtheorem{theorem}{Theorem}

\newtheorem{corollary}{Corollary}
\newtheorem{proposition}{Proposition}
\newtheorem{definition}{Definition}

\def \Xr {\mathcal{X}_r}
\def \P{\mathfrak{P}}
\def \Proj{\mathbb{P}}
\def \N {\ensuremath{\mathbb{N}}}

\def \mcF {\ensuremath{\mathcal{F}}}
\def \oK{\mathfrak{o}_K}

\def \P {\mathfrak{P}}
\def \Q {\mathbb{Q}}

\def \PG {\ensuremath{\mathfrak{P}_G}}

\def \isom {\ensuremath{\stackrel{\sim}{=}}}

\def \T {\ensuremath{\mathfrak{T} }}
\def \Ttwo {\ensuremath{\mathfrak{T}_2 }}

\def \Qp {\mathbb{Q}_p}
\def \Z {\mathbb{Z}}
\def \z {\zeta^\triangleleft}
\def \Zp  {\mathbb{Z}_p}

\def \Fp  {\mathbb{F}_p}
\def \Fq  {\mathbb{F}_q}

\def \F23 {\ensuremath{F_{2,3}}}

\def \bfg { \overline{{\bf g}}}
\begin{document}
\maketitle

\begin{abstract}
We give explicit formulae for the local normal zeta functions of
torsion-free, class-$2$-nilpotent groups, subject to conditions on
the associated Pfaffian hypersurface which are generically satisfied
by groups with small centre and sufficiently large abelianization. We show how the
functional equations of two types of zeta functions - the
Weil zeta function associated to an algebraic variety and 
zeta functions of algebraic groups introduced by Igusa - match up to
give a functional equation for local normal zeta functions of groups. We also give explicit formulae and derive functional equations for an
 infinite family of class-$2$-nilpotent groups known as Grenham
 groups,
confirming conjectures of du~Sautoy.\footnote{2000 AMS Mathematics Subject Classification 11 M 41, 20 E 07.}
\end{abstract}

\section{Introduction}
In~\cite{GSS/88} Grunewald, Segal and Smith introduced the concept of
a zeta function of an infinite, finitely generated, torsion-free nilpotent group~$G$
(a \T-group, in short). To any family~$\mathcal{X}$ of subgroups of~$G$ they
associate the abstract Dirichlet series
$$\zeta_{\mathcal{X}}(s)=\sum_{H \in
  \mathcal{X}}|G:H|^{-s}=\sum_{n=1}^\infty a_n(\mathcal{X})n^{-s},$$
where 
$$ a_n(\mathcal{X})=|\{H\in\mathcal{X}|\; |G:H|=n\}|$$
(and $\infty^{-s}=0$). Among the families whose study was initiated in~\cite{GSS/88} are 
\begin{eqnarray*}
\mathcal{S}(G)&=&\{\text{all subgroups of finite index in }G\},\\
\mathcal{N}(G)&=&\{\text{all normal subgroups of finite index in }G\},\\
\mathcal{H}(G)&=&\{H\in\mathcal{S}(G)|\; \widehat{H}\isom \widehat{G}\},
\end{eqnarray*}
where $\widehat{G}$ denotes the profinite completion of the
group~$G$. These give rise to the zeta functions
\begin{eqnarray*}
\zeta^\leq_G(s)&=&\zeta_{\mathcal{S}(G)}(s),\\
\zeta^\triangleleft_G(s)&=&\zeta_{\mathcal{N}(G)}(s),\\
\zeta^{\widehat{\;}}_G(s)&=&\zeta_{\mathcal{H}(G)}(s).
\end{eqnarray*}
All these zeta functions satisfy an Euler product decomposition
$$ \zeta_{G}^*(s)=\prod_{p \text{ prime}} \zeta_{G,p}^*(s),\;*\in\{\leq,\triangleleft,\widehat{\;}\},$$
into a product of {\sl local} zeta functions, were
$$\zeta^*_{G,p}(s)=\sum_{n=0}^\infty a^*_{p^n}p^{-ns}$$
and
$$a^\leq_{p^n}=a_{p^n}(\mathcal{S}),\;a^\triangleleft_{p^n}=a_{p^n}(\mathcal{N}),\;a^\wedge_{p^n}=a_{p^n}(\mathcal{H}).$$
 It is the local
  normal zeta functions~$\zeta_{G,p}^\triangleleft(s)$ we will concentrate on in
this paper. 

One of the main results
of~\cite{GSS/88} establishes the rationality of the local zeta
functions in~$p^{-s}$. It relies on the presentation of local zeta
functions as $p$-adic integrals:

\begin{theorem}[\cite{GSS/88}, Prop. 3.1. and Thm. 4.1] Given a $\T$-group~$G$ of
  Hirsch length\footnote{ Recall that the {\sl Hirsch length} of a~$\T$-group~$G$ is
the number of infinite cyclic factors in a decomposition series
for~$G$. } $h(G)=n$. For
  $*\in\{\leq,\triangleleft,\widehat{\;}\}$ and almost all primes~$p$
$$\zeta_{G,p}^*(s)=(1-p^{-s})^{-n}\int_{V^*_p}|m_{11}|^{s-1}\dots|m_{nn}|^{s-n}|dx|$$
where $|m|=p^{-v_p(m)}$, $v_p$ is the valuation on~$\Zp$, $|dx|$ is
the normalized additive Haar measure on~$\Zp^{d(d+1)/2}\equiv
Tr_n(\Zp),$ the triangular $n\times n$-matrices over~$\Zp$ with
diagonal entries~$m_{ii}$ and suitable subsets~$V^*_p\subseteq Tr_n(\Zp)$.
\end{theorem}
Rationality is a consequence of the observation that the
subsets~$V^*_p$ are ${\sl definable}$ in the language of fields. In
this situation, a theorem of Denef's~\cite{Denef/84} is applicable which in turn
relies on an application of Macintyre's quantifier elimination for the
theory of~$\Qp$~\cite{Macintyre/76} and on Hironaka's theorem~\cite{Hironaka/64} on
resolution of singularities in characteristic zero. 

A major challenge in the field is to understand how the local
zeta functions vary with the prime~$p$. The zeta
function~$\zeta_{G}^*(s)$ is called {\sl finitely uniform} if there
are finitely many rational functions~$V_i(X,Y)\in\mathbb{Q}(X,Y)$,
$1\leq i\leq r$, such that for each prime~$p$ there exists an $i$ such
that 
\begin{equation}
\zeta_{G,p}^*(s)=V_i(p,p^{-s}),\label{fin-uniformity}
\end{equation}
and {\sl uniform} if~$r=1$. 

(Finite) uniformity is not typical for zeta functions of nilpotent groups,
however: Du Sautoy and Grunewald linked the question of the local
factors' dependence on the prime~$p$ to the classical problem of counting points on varieties
mod~$p$. In~\cite{duSG/00} they identify local zeta functions of
groups as special cases of a more general class of $p$-adic integrals
they call {\sl cone integrals}:
\begin{definition}[\cite{duSG/00}, Def. 1.2] Let $f_i,g_i$,
  $i\in\{0,1,\dots,l\}$ be polynomials over~$\Q$ in~$n$ variables. The condition
\begin{equation}
\psi({\bf x}):\;v_p(f_i({\bf x}))\leq v_p(g_i({\bf x})) \text{ for
}i\in\{1,\dots,l\}\label{cone condition}
\end{equation}
is called a {\rm cone condition}. A $p$-adic integral of the form
$$Z_{({\bf f},{\bf g})}(s,p)=\int_{\{{\bf x}\in\Zp^n|\;\psi({\bf x})
  \text{ holds}\}}|f_0({\bf x})|^{s}|g_0({\bf x})||dx|$$
is called a {\rm cone integral}. The vector $({\bf f},{\bf g})$ is called {\rm
  cone integral data}.
\end{definition}
If the condition~(\ref{cone condition}) is trivial (and~$g_0\equiv1$) we recover Igusa's
local zeta function (cf. Appendix to~\cite{Igusa/77}, \cite{Denef/87})
as a special case of a cone integral.
Writing local zeta functions of nilpotent groups as cone integrals not
only allowed the authors of~\cite{duSG/00} to dispense with the -~in
general mysterious~- model-theoretic black box of quantifier
elimination. It also enabled them to give an -~in principal~- very
explicit expression for these functions.
\begin{theorem}\label{duSG theorem} [Thm 1.6 in~\cite{duSG/00}] Let~$G$ be a~\T- group,
  and $*\in\{\leq,\triangleleft,\widehat{\;}\}.$ There exists an algebraic
  variety~$Y^*$ defined over~$\Q$, with irreducible
  components~$E^*_i,\,i\in T^*:=\{1,\dots,t^*\}$, all of which are smooth and intersect
  normally, and rational
functions~$P^*_I(X,Y)\in\mathbb{Q}(X,Y)$, $I\subseteq T^*$ such that, for almost all primes~$p$,
\begin{equation}
 \zeta^*_{G,p}(s)=\sum_{I\subseteq
  T^*}c^*_{p,I}P^*_{I}(p,p^{-s})\label{duSG formula}
\end{equation}
where
$$
c^*_{p,I}=|\{a\in\overline{Y^*}(\Fp):\;a\in\overline{E^*_i}\mbox{
  if and only if }i\in I\}|$$
and $\overline{Y}$ means the reduction$\mod p$ of the variety~$Y$.
\end{theorem}
The varieties~$Y^*$ arise as resolutions of singularities of the
hypersurfaces~$(\prod_{i=0}^lf^*_ig^*_i=0)$, where~$({\bf f}^*,{\bf g}^*)$ is
the cone integral data for the respective cone
integrals. Theorem~\ref{duSG theorem} suggests that finite uniformity
should be the exception rather than the rule for zeta functions of nilpotent groups. But the question
what varieties may appear in~(\ref{duSG formula}) in general remains wide open. Only recently du Sautoy presented the first example of a \T-group~$G$ for
which neither~$\zeta_{G}^\leq(s)$ nor $\zeta_{G}^\triangleleft(s)$ are
finitely uniform, but depend on the number of~$\Fp$-points of an
elliptic curve~(\cite{duS-ecI/01},\cite{duS-ecII/01}).

The first explicit formulae for 
  non-uniform\footnote{By abuse of language we will say
  `non-uniform' for `not finitely uniform'.}
  normal zeta functions, including du Sautoy's examples, appeared
  in~\cite{Voll/03}. Rather than evaluating cone integrals we
  introduced in this paper a much less coordinate-dependent calculus
  to compute local normal zeta functions of $\T$-groups of nilpotency
  class~$2$. It is based on an enumeration of
  vertices in the affine Bruhat-Tits buildings associated to~$Sl_n(\Qp)$. In
  Section~\ref{section-lattices-and-flags} we will recall
  from~\cite{Voll/03} what is necessary to make the current paper
  self-contained.

The examples computed in~\cite{Voll/03} suggested that non-uniform
  zeta functions, too, might satisfy certain {\sl local functional
  equations}. The (finitely uniform) examples computed before had all
  featured a functional equation of the form
\begin{equation}
V_i(X^{-1},Y^{-1})=(-1)^{l_i}\, X^{m_i}Y^{n_i}\,V_i(X,Y)\label{combinatorial-funeq}
\end{equation}
for integers $l_i,m_i,n_i$ and $V_i(X,Y)$ as in~(\ref{fin-uniformity})
above, which defied explanation for~$*\in\{\leq,\triangleleft\}$. In the non-uniform examples given
in~\cite{Voll/03}, the uniform components showed symmetries
like~(\ref{combinatorial-funeq}), too, matching with the functional
equation of the Weil zeta function counting the number
of~$\mathbb{F}_q$-points of algebraic varieties to give a functional
equation for the local zeta functions. But the nature of these uniform
components and their functional equations remained poorly understood. In this paper we try to shed some light
onto this phenomenon. 

\section{Statement of results}
Let now $G$ be a \T-group of nilpotency class~$2$ (a \Ttwo-group, in short) with
derived group~$G':=[G,G]$ and centre~$Z(G)$. Only
for simplicity we make the following
\begin{assumption}$G/G'$ and $G'$ are {\rm  torsion-free} abelian of rank~$d$
  and~$d'$, respectively, and~$G=Z(G)$. \label{assumption}
 \end{assumption}
Indeed, in general both
$Z(G)/G'$ and $G/G'$ are finitely generated abelian groups. But as we are only
looking to prove results about all but finitely many of the local zeta functions
$\zeta^\triangleleft_{G,p}(s)$ we loose nothing by restricting
ourselves to primes~$p$ not dividing the orders of the
respective torsion parts. And as one checks with no difficulty that
 $$\zeta^\triangleleft_{G\times
      \Z^r}(s)=\zeta^\triangleleft_{G}(s)\cdot\zeta(s-n)\zeta(s-(n+1))\dots\zeta(s-(n+r-1)),$$ where $h(G)=n$ and $\zeta(s)=\sum_{k=1}^\infty k^{-s}$ is the Riemann zeta function, we may indeed assume that~$Z(G)=G'$. 
Under Assumption~\ref{assumption}~ the group~$G$ has a presentation
\begin{equation}
 G=\langle x_1,\dots,x_d,y_1,\dots,y_{d'}|\;[x_i,x_j]=M({\bf
  y})_{ij}\rangle, \label{presentation}
\end{equation}
where $M({\bf y})$ 
 is an anti-symmetric $d\times d$-matrix of $\mathbb{Z}$-linear
forms in ${\bf y}=(y_1,\dots,y_{d'})$, and all other commutators are
trivial. Note that we adopted additive notation for words
in~$Z(G)$. Conversely, of course, every such matrix~$M({\bf y})$
defines a $\Ttwo$-group via~(\ref{presentation}). 

If the polynomial $\text{Pf}\,(M({\bf y})):=\sqrt{\det(M({\bf
    y}))}$ is not identical zero we call the hypersurface~$\PG$
in~$\mathbb{P}^{d'-1}$ defined over the integers by $\text{Pf}\,(M({\bf y}))=0$ the {\sl Pfaffian hypersurface} associated
to~$G$. Given a fixed prime number~$p$ we then denote by~$\overline{\PG}$ its reduction modulo~$p$. Our main result is

\begin{theorem} \label{theorem-funeq} Assume that ${\rm Pf}(M({\bf
  y}))\in\Z[{\bf y}]$ is non-zero and irreducible. Assume that
  the Pfaffian hypersurface~$\PG$ is smooth and contains no lines. For
  a prime~$p$ let $$n_{\PG}(p)=|\PG(\Fp)|$$ denote the number of
  $\Fp$-rational points of~$\PG$. Then there are (explicitly
  determined) rational functions $W_{0}(X,Y),
  W_{1}(X,Y)\in\mathbb{Q}(X,Y)$ such that if~$\PG$ has good reduction mod~$p$
\begin{equation}
\zeta^\triangleleft_{G,p}(s)=W_{0}(p,p^{-s})+n_{\PG}(p)\cdot
W_{1}(p,p^{-s}). \label{theorem-funeq-expression}
\end{equation}\end{theorem}

\begin{corollary}\label{coro to main theorem}
If, moreover,~${\PG}$ is  absolutely irreducible, the following functional equation holds:
\begin{equation}
\zeta^\triangleleft_{G,p}(s)|_{p\rightarrow p^{-1}}=(-1)^{d+d'}p^{\binom{d+d'}{2}-(2d+d')s}\zeta^\triangleleft_{G,p}(s).\label{funeq}
\end{equation}
\end{corollary}

We should like to remark that the conditions of
  Theorem~\ref{theorem-funeq} on the $\Ttwo$-group~$G$ are
  {\sl generically} satisfied for small~$d'$ and large $d$\;($=2r$,
  say):  A presentation as in~(\ref{presentation}) specifies a $\Z$-linear embedding of a $\Proj^{d'-1}$ into
  the projective space $\Proj(\mathcal{S})$ over the vector space~$\mathcal{S}$ of anti-symmetric $(d\times d)$-matrices. The Pfaffian
  hypersurface~$\PG$ is just the intersection of this~$\Proj^{d'-1}$
  with the universal Pfaffian
  hypersurface~$\mathcal{X}_r\subset\Proj(\mathcal{S})$ of singular
  matrices. The singular locus of the latter consists of matrices of rank~$\leq 2r-4$ and has
  codimension~$6$. A {\sl generic} $\Proj^{d'-1}\subset\Proj(\mathcal{S})$ therefore intersects~$\mathcal{X}_r$ along a {\sl smooth}
  hypersurface in~$\Proj^{d'-1}$ of degree~$r$ if $d'\leq 6$ (This is
  remark~(8.3) in~\cite{Beauville/00}). 

More work is required to see that for all~$d'$ and
  for $d>4d'-10$ a generic Pfaffian hypersurface~$\PG$ will not
  contain lines. This follows immediately from the following
  proposition. It is due to
  Arnaud Beauville and we cordially thank him for his contribution of its proof as an Appendix to this paper.

\begin{proposition}\label{beauville-proposition} A generic {\rm Pfaffian} hypersurface of degree~$r>2n-3$ in~$\mathbb{P}^n$
  contains no lines.
\end{proposition}

In view of these remarks it seems to us as if {\sl smoothness} of the
Pfaffian~$\PG$ was the
most restrictive condition on the~\Ttwo-group~$G$ in
Theorem~\ref{theorem-funeq}. In~\cite{Voll/03b} we examine examples
showing that for singular Pfaffians~$\PG$ the expressions for the zeta
functions of~$G$ will reflect the rank-stratification of the
determinantal variety~$\PG$. In
Section~\ref{section-lattices-and-flags} we shall comment briefly on
the technical reasons for the condition on the Pfaffian hypersurface
to contain no lines. The first inroads in
overcoming these obstacles were made by Pirita Paajanen, a student of
du~Sautoy, who derived the following explicit formula for the local normal zeta functions of~$F_{2,4}$, the
free class-$2$-nilpotent group on four generators. In this example the Pfaffian
hypersurface~$\PG$ is a smooth quadric fourfold.

\begin{proposition}\label{Paajanen proposition}   [Paajanen, \cite{PaajanenDPhil}] Let $G=F_{2,4}$. Then for all primes $p$
$$ \zeta_{G,p}^\triangleleft(s)=W_0(p,p^{-s})
                              +n^{(1)}_{\PG}(p)W_1(p,p^{-s})
                              + n^{(2)}_{\PG}(p)W_2(p,p^{-s})
                              + n^{(3)}_{\PG}(p)W_3(p,p^{-s}),$$
where 
\begin{eqnarray*}
n^{(1)}_{\PG}(p)&=&(p^2+1)(p^2+p+1),\\
      n^{(2)}_{\PG}(p)&=&(p+1)(p^2+1)(p^2+p+1),\\
      n^{(3)}_{\PG}(p)&=&2(p^2+1)(p+1)
\end{eqnarray*}
denotes the number of $\Fp$-rational points of
the Fano varieties of $(i-1)$-dimensional subspaces on~$\PG$, and
$W_i(X,Y)$ are rational functions. The functional
equation~(\ref{funeq}) holds.
\end{proposition}

The proof of Theorem~\ref{theorem-funeq} will be given in two steps:
First we prove it in the
case $n_{\PG}(p)=0$ to obtain the rational
function~$W_{0}(p,p^{-s})$. Note that there are \Ttwo-groups for which
$n_{\PG}(p)=0$
  for infinitely many primes~$p$. (Indeed, if~ $G=H(\oK)$ is the Heisenberg group over the ring of integers of an
  algebraic number field~$K$, these are exactly the {\sl inert}
  primes. The remaining primes are harder to deal with as the
  associated Pfaffian will not be smooth.) Loosely speaking, if~$p$ is a prime for which~$\PG$ defines a smooth non-empty
hypersurface in projective~$(d'-1)$-space over~$\Fp$ without lines we will need to
`correct' the rational function $W_{0}(p,p^{-s})$ along
the~$n_{\PG}(p)$ $\Fp$-points of~$\PG$ by the
function~$W_{1}(p,p^{-s})$ to obtain the $p$-th local normal zeta function. This will constitute
the second step in the proof of Theorem~\ref{theorem-funeq}.

To prove Corollary~\ref{coro to main theorem} we shall demonstrate
that the functional equation~(\ref{funeq}) is due to the interplay of
two phenomena: We will firstly recall how the {\sl functional equation} of the
Hasse-Weil zeta function associated to the hypersurface~$\PG$ gives
rise to a symmetry of the expression~$n_{\PG}(p)$ as a function of~$p$
(and how this zeta function's {\sl rationality} gives sense to the
symbol~$n_{\PG}(p^{-1})$ in the left hand side of~(\ref{funeq})in the first place).  
Secondly we shall show that the uniform components~$W_i(X,Y)$,
$i\in\{0,1\}$, in~(\ref{theorem-funeq-expression})
satisfy a symmetry of the form~(\ref{combinatorial-funeq}). Our main
tool to achieve this will be
Theorem~\ref{theorem-flag-rational-function} below, which is
essentially due to Igusa~\cite{Igusa/89}. It establishes
such a symmetry for a single rational function
(\ref{formula-flag-rational-function}) which is defined in terms of flag varieties.

\begin{theorem} \label{theorem-flag-rational-function}
Let $n\geq 2$ be an integer, $X_1,\dots,X_{n-1}$ independent indeterminates
and~$q$ a prime power. For $I\subseteq\{1,\dots,n-1\}$ let $b_I(q)=|\mcF_I(\Fq)|\in\mathbb{Z}[q]$ denote the number of $\Fq$-points of the projective variety of flags in $\Fq^n$ of type~$I$. Set
\begin{equation}\label{formula-flag-rational-function}
F_n(q,{\bf X}):=\sum_{I\subseteq\{1,\dots,n-1\}}b_I(q)\prod_{i\in I}\frac{X_i}{1-X_i}.
\end{equation}
Then
\begin{equation}
F_n(q^{-1},{\bf X}^{-1})=(-1)^{n-1}q^{-\binom{n}{2}}F_n(q,{\bf X}).\label{funeq-rat-function}
\end{equation}
\end{theorem}
(See Section~\ref{section-lattices-and-flags} for the definition of
`flag of type~$I$'.)

We will show how (the crucial factors of) the rational
functions~$W_i(p,p^{-s})$ in~(\ref{theorem-funeq-expression})
may be derived from functions of type~(\ref{formula-flag-rational-function}) by
suitable substitutions of variables. The~$W_i(p,p^{-s})$ `inherit' the
functional equation~(\ref{funeq-rat-function}) as the quotient
$F_n(q^{-1},{\bf X}^{-1})/F_n(q,{\bf X})$  is independent of the
`numerical data'~${\bf X}$.

The zeta functions introduced by Igusa~\cite{Igusa/89} are
defined in terms of root systems of algebraic groups. In fact we only
need the most basic of these, the one associated to~$Gl_n$. Our
formulation (and our elementary proof) of
Theorem~\ref{theorem-flag-rational-function} in the language of flag
varieties seems natural from the point of view taken
in~\cite{duSG/00}, whereas a connection to algebraic groups seems
elusive in the context of normal zeta functions of groups.

In~\cite{duSLubotzky/96}, however, du Sautoy and
Lubotzky interpret the zeta
functions $\zeta^{\widehat{\;}}(s)$ (where~$G$ is again a general \T-group) as~$p$-adic integrals over the algebraic automorphism group of the Lie algebra
associated to~$G$. A generalisation of Igusa's work~\cite{Igusa/89}
allows them to derive uniformity as well as local functional equations of
these zeta functions for certain classes of \T-groups.

In~\cite{DenefMeuser/91} Denef and Meuser prove a functional equation for the Igusa
local zeta function associated to a homogeneous polynomial. These zeta functions also have an expression
of the form~(\ref{duSG formula}). However, the `uniform
components' occurring in this context show rather less
structure than the rational functions~$W_i(X,Y)$ in Theorem~\ref{theorem-funeq} of the current paper.

\bigskip

As another application of Theorem~\ref{theorem-flag-rational-function}
we derive both  explicit formulae and  local functional equations for
the normal zeta functions of another infinite family
of~$\Ttwo$-groups known as `Grenham's groups' (cf. \cite{duS/00},
Chapter~5.8, or~\cite{duSG-ghosts/00}, Chapter~6.3). For~$n\geq2$ let 
$$G_n:=\langle  x_1,\dots,x_n,y_1,\dots,y_{n-1}|\;[x_i,x_n]=y_i, 1\leq
i\leq n, \text{ all other }[,]\text{ trivial}\rangle.$$
The group $G_n$ may be thought of as~$n-1$ copies of the discrete
Heisenberg group -~the group of~$3\times3$-upper uni-triangular
matrices with integer entries~- with one off-diagonal entry identified
in each copy.
\begin{theorem}\label{theorem-grenham} For all primes~$p$
$$\z_{G_n,p}(s)=\zeta_{\Zp^n}(s)\zeta_p((2n-1)s-n(n-1))F_{n-1}(p^{-1},{\bf X})$$
where ${\bf X}=(X_1,\dots,X_{n-2})$ and
$$X_i=p^{-(2(n-i)-1)s+(n+i)(n-i-1)}\text{ for }i\in\{1,\dots,n-2\}.$$
In particular, the following functional equation holds:
\begin{equation}
\zeta^\triangleleft_{G_n,p}(s)|_{p\rightarrow
  p^{-1}}=-p^{\binom{2n-1}{2}-(3n-1)s}\zeta^\triangleleft_{G_n,p}(s).\label{grenham
  normal funeq}
\end{equation}
\end{theorem}
In the forthcoming paper~\cite{Voll/04} we use our method to compute explicitly
all the {\sl subgroup} zeta functions~$\zeta^\leq_{G_n}(s)$ and prove that
\begin{equation}
\zeta_{G_n,p}^{\leq}(s)|_{p\rightarrow
  p^{-1}}=-p^{\binom{2n-1}{2}-(2n-1)s}\zeta_{G_n,p}^{\leq}(s).\label{grenham  subgroup funeq}
\end{equation}
The functional equations~(\ref{grenham normal funeq})
and~(\ref{grenham subgroup funeq}) were conjectured by
du~Sautoy (Conjecture~5.41 in~\cite{duS/00}).

We will prove Theorem~\ref{theorem-flag-rational-function} in Section~\ref{subsection theorem rat function} by an argument using the
Schubert cell decomposition of flag varieties. The proofs of
Theorem~\ref{theorem-funeq} and
Theorem~\ref{theorem-grenham} will occupy Sections~\ref{subsection proof of theorem-funeq} and~\ref{subsection grenham}, respectively. The point of view taken is the one
developed in~\cite{Voll/03} (where also the special cases of
Theorem~\ref{theorem-funeq} for~$d'\in\{2,3\}$ were proved). There the {\sl Cartan decomposition} for
lattices in the centre of~$G$ was used to interpret the local zeta
functions as generating functions associated to certain weight
functions on the vertices of the Bruhat-Tits building~$\Delta_{d'}$
for~$Sl_{d'}(\mathbb{Q}_p)$, exhibiting its dependence on the geometry
of the Pfaffian hypersurface~$\PG$. We will recall briefly the main results
of~\cite{Voll/03} in Section~\ref{section-lattices-and-flags} together
with some basic definitions and observations about lattices and flags.

\begin{acknowledgements} We should like to thank Konstanze
  Rietsch, Fritz Grunewald and Marcus du Sautoy
  for helpful and inspiring discussions. The suggestions made by a
  referee were a great help in improving the exposition of this paper. We
  gratefully acknowledge support by the UK's Engineering and Physical Sciences Research Council (EPSRC) in form of a Postdoctoral Fellowship. 
\end{acknowledgements}

\section{Flags and lattices}\label{section-lattices-and-flags}
In this section we give a brief summary of the method developed
in~\cite{Voll/03} to compute local normal zeta functions.

Let $G$ be a $\Ttwo$-group satisfying Assumption~\ref{assumption}.
For a fixed prime~$p$, the computation of the $p$-th normal local zeta
function of~$G$ comes down to an enumeration of those lattices in the $\Zp$ -Lie algebra (with Lie brackets induced by taking commutators) $$G_p:=(G/Z(G)\oplus Z(G))\otimes_{\Z} \Zp$$ which
 are {\sl ideals} in $G_p$. We
call a lattice $\Lambda\subseteq\Zp^n$ {\sl
  maximal} (in its homothety class) if $p^{-1}\Lambda\not\subseteq\Zp^n$. The key observation is the following
\begin{lemma}\label{gss-lemma}[Lemma
6.1 in~\cite{GSS/88}]
For each lattice $\Lambda'\subseteq G_p'$ put $X(\Lambda')/\Lambda'=Z(G_p/\Lambda')$. Then
\begin{eqnarray*}
\lefteqn{\zeta^\triangleleft_{G,p}(s)=\zeta_{\Zp^{d}}(s)\sum_{\Lambda'\subseteq
  G_p'}|G_p':\Lambda'|^{d-s}|G_p:X(\Lambda')|^{-s}}\\
&=&\zeta_{\Zp^{d}}(s)\zeta_p((d+d')s-dd')\underbrace{\sum_{\substack{\Lambda'\subseteq G'_p  \\ \Lambda' {\text maximal}}}|G_p':\Lambda'|^{d-s}|G_p:X(\Lambda')|^{-s}}_{=:A(p,p^{-s}),{\text say}}
\end{eqnarray*}
\end{lemma}
\begin{corollary}\label{coro to gss-lemma}
\begin{eqnarray*}
\left.A(p,p^{-s})\right|_{\substack{p \rightarrow p^{-1}}}&=&(-1)^{d'-1}p^{\binom{d'}{2}} A(p,p^{-s}) \label{fun.eq.buildings}\\
\quad\Longleftrightarrow\quad\z_{G,p}(s)|_{p\rightarrow
  p^{-1}}&=&(-1)^{d+d'}p^{\binom{d+d'}{2}-(2d+d')s}\z_{G,p}(s).
\label{fun.eq.} 
\end{eqnarray*}\hfill $\Box$
\end{corollary}
 Recall that (homothety\footnote{Two lattices
   $\Lambda,\Lambda'\subseteq \Qp^{d'}$ are called {\sl homothetic} if
   there is a non-zero constant $c\in\Qp$ such that
   $c\Lambda=\Lambda'$.} classes defined by) maximal lattices are in one-to-one
correspondence with vertices of the Bruhat-Tits building $\Delta_{d'}$
for $Sl_{d'}(\Qp)$ (e.g.~\cite{Garrett/97}, \S 19).  To derive an
explicit formula for $A(p,p^{-s})$ requires a quantitative
understanding of two integer-valued functions~$w$ and~$w'$ on the vertex set of the
simplicial complex~$\Delta_{d'}$. We write
\begin{equation}
A(p,p^{-s})=\sum_{[\Lambda']}p^{dw([\Lambda'])-sw'([\Lambda'])},\label{A}
\end{equation}
where, for a homothety class $[\Lambda']$ of a {\sl maximal} lattice
$\Lambda'$ in $G_p'\isom \Zp^{d'}$ we define
\begin{eqnarray*}
w([\Lambda']) &:=& \log_p(|G'_p:\Lambda'|),\\
w'([\Lambda']) &:=& w([\Lambda'])+ \log_p(|G_p:X(\Lambda')|).
\end{eqnarray*}
In order to describe the dependence of these functions on the
 lattice~$\Lambda'$ we will introduce some notation. For an integer~$m\in\mathbb{N}_{>0}$ we set~$[{m}]:=\{1,\dots,m\}$. A
 maximal lattice~$\Lambda'\subseteq\Zp^{d'}$ is said to be of {\sl
   type}~$\nu(\Lambda')=(I,{\bf r}_I)$ if 
\begin{equation}
I=\{i_1,\dots,i_{l}\}\subseteq[{d'-1}],\;{\bf r}_I=(r_{i_1},\dots,r_{i_{l}}
),\;i_1<\dots <i_l,\label{I}
\end{equation}
and $r_{i_j}\in\mathbb{N}_{>0}\mbox{ for }j\in\{1,\dots,l\}$ and $\Lambda'$ has elementary divisors
\begin{gather}
\left(\underbrace{1,\dots,1}_{i_1},\underbrace{p^{r_{i_1}},\dots,p^{r_{i_1}}}_{i_2-i_1},\underbrace{p^{r_{i_1}+r_{i_2}},\dots,p^{r_{i_1}+r_{i_2}}}_{i_3-i_2},\dots,\underbrace{p^{\sum_{j=1}^{l}r_{i_j}},\dots,p^{\sum_{j=1}^{l}r_{i_j}}}_{d'-i_{l}}\right)=:(p^{\nu}).\label{eldiv type}
\end{gather}
By slight abuse of notation we may  say that a maximal lattice is of
type~$I$ if it is of type~$(I,{\bf r}_I)$ for some positive
vector~${\bf r}_I$ and that the {\sl homothety class}~$[\Lambda']$ has
type~$I$ if its maximal element has type~$I$, in which case we
write~$\nu([\Lambda'])=I$. For computations it shall be advantageous
to write
\begin{eqnarray}
A(p,p^{-s})&=&\sum_{I\subseteq[d'-1]}A_I(p,p^{-s}),\text{ where}\label{A_d'-sum-over-types}\\
A_I(p,p^{-s})&:=&\sum_{\nu([\Lambda'])=I}p^{dw([\Lambda'])-sw'([\Lambda'])}.\nonumber
\end{eqnarray}
Notice that the lattice's index -~and thus~$w([\Lambda'])$~- is given by
\begin{equation}
|\Zp^{d'}:\Lambda'|=p^{\sum_{i\in
 I}r_{i}(d'-i)}=p^{w([\Lambda'])}\;.\label{index}
\end{equation}

We shall explain how the evaluation of~$w'$ may be
reduced to solving linear congruences. 
The group $\Gamma=Sl_{d'}(\Zp)$ acts transitively on the set of maximal lattices of fixed type. If
we choose a basis for the $\Zp$-module $G_p'$, represent lattices
as the row span of $d'\times d'$-matrices and denote by $\Gamma_{\nu}$
the stabilizer in $\Gamma$ of the lattice generated by the diagonal
matrix whose entries are given by the vector~(\ref{eldiv type}), the orbit-stabiliser
theorem gives us a $1-1$ correspondence
\begin{equation}
\left\{\mbox{maximal lattices of type }(I,{\bf
    r}_I)\right\}\stackrel{1-1}{\longleftrightarrow}
    \Gamma/{\Gamma_{\nu}}.\label{lattice correspondence}
\end{equation}
This correspondence allows us to describe 
$|G_p:X(\Lambda')|$ -~and thus~$w'([\Lambda'])$~- for a maximal lattice~$\Lambda'$ in terms of $M({\bf
  y})$, the matrix of commutators in a presentation for~$G$ as
in~(\ref{presentation}). 

\begin{theorem}\label{theorem linear congruences}[\cite{Voll/03}, \S 2.2] Let~$\Lambda'$ correspond
  to the coset~$\alpha\Gamma_{\nu}$ under~(\ref{lattice correspondence}), where $\alpha\in\Gamma$ with
{\rm column} vectors $\alpha^j$, $j=1,\dots,d'$. Then 
$|G_p:X(\Lambda')|$ equals the index of the kernel of the following system
of linear congruences in~$G_p/G_p'$:
\begin{equation}
\forall i \in \{1,\dots,d'\}\quad {\bf \overline{g}}M(\alpha^i)\equiv
0 \; mod \;(p^{\nu})_i,\label{linear congruences}
\end{equation}
where $ {\bf \overline{g}}=(\overline{g}_1,\dots,\overline{g}_d)\in
G_p/G'_p\isom \Zp^d$ and $(p^{\nu})_i$ denotes the $i$-th entry of the
vector $(p^{\nu})$ given in~(\ref{eldiv type}).
\end{theorem}
A {\sl flag of type $I$} in~$\mathbb{P}^{d'-1}(\Fp)$,
$I\subseteq[{d'-1}]$ as in~(\ref{I}), is a sequence~$(V_i)_{i\in I}$ of incident vector spaces
$$\mathbb{P}^{d'-1}(\Fp)>V_{i_1}>\dots> V_{i_{l}}>\{0\}$$
with~$\mbox{{\bf co}dim}_{\Fp}(V_i)=i$.
A flag is called {\sl incomplete} (or {\sl partial}) if $I\not=[{d'-1}]$, and {\sl complete} otherwise.
The flags of type $I$ form a projective variety $\mcF_I$, whose number
of $\Fp$-points is given by $b_I(p)\in\mathbb{Z}[p]$, a polynomial
whose leading term equals~$p^{\dim{\mcF_I}}$. These polynomials are
easily expressed in terms of~$p$-binomial coefficients, but we will
not make use of this fact here. It is easy to see that for all $I\subseteq
[{n-1}]$

\begin{equation}
b_I(p^{-1})=p^{-\dim{\mcF_I}}b_I(p).\label{funeq-varieties}
\end{equation}

\noindent Given~$\nu=(I,{\bf r}_I)$ as above, let $f(I,{\bf r}_I,p)=|\Gamma/\Gamma_{\nu}|$
be the number of maximal lattices in~$\Zp^{d'}$ of type~$\nu$. Using~(\ref{lattice correspondence}) one proves easily
\begin{lemma}\label{lemma-lattices-of-given-type}
$$f(I,{\bf r}_I,p)=b_I(p)\; p^{\sum_Ir_{i}(d'-i)i-\dim\mcF_I}.$$
\end{lemma}
 Given $\alpha\in\Gamma$ and~$I$ as in~(\ref{I}), let
$\overline{\alpha}$ denote the reduction mod~$p$ and define vector spaces
$$V_{i}:=\langle
\overline{\alpha^{i+1}},\dots,\overline{\alpha^{d'}}\rangle < \mathbb{P}^{d'-1}(\Fp)  ,\quad i\in I.$$
Clearly ${\text codim}(V_{i})=i$. We will call the
flag~$(V_{i})_{i\in I}$ of type~$I$ {\sl the flag associated to} $\Lambda'$ if $\nu(\Lambda')= I $ and $\Lambda'$
corresponds to $\alpha\Gamma_{\nu}$
under~(\ref{lattice correspondence})\footnote{It is indeed straightforward to show that this
    is well-defined, i.e. independent of the coset representative
    $\alpha$.}. 
Given a fixed point ${x}\in\mathbb{P}^{d'-1}(\Fp)$, we call a
maximal lattice~$\Lambda'$ {\sl a lift of}~${x}$ if its associated
flag contains~${x}$ as $0$-dimensional member and we shall
write~${x}=P([\Lambda'])$. Note that then
necessarily $d'-1\in\nu([\Lambda'])$.

We can now explain why we made the assumption that~$\PG$ should be
smooth and contain no lines. In the latter case,~(\ref{linear congruences}) is
equivalent to
\begin{eqnarray}\bfg&\equiv& 0 \mod p^{\sum_{i\in
      I\setminus\{d'-1\}}r_i}\nonumber\\
 p^{-\sum_{i\in
      I\setminus\{d'-1\}}r_i}\bfg M(\alpha^{d'})&\equiv&0\mod
 p^{r_{d'-1}}\label{equiv congruence}
\end{eqnarray}
(where the congruence~(\ref{equiv congruence}) is regarded trivial if~$d'-1\not\in
I$). This follows easily from the observation that under this
assumption the vectors $\alpha^1,\dots,\alpha^{d'-1}$ may always be
chosen to lie outside the Pfaffian hypersurface, i.e. such that
$v_p(\det(M(\alpha^i)))=0$ for $i\in\{1,\dots,d'-1\}$. In the case
$d'-1\in I$ we are left to analyse the elementary divisors of the
matrix~$M(\alpha^{d'})$. But by an easy geometrical argument
(cf. Lemma~1, \cite{Voll/03b}) one sees that if the Pfaffian is smooth
this matrix always has a
$(d-2)\times (d-2)$-unit minor\footnote{This may or may not be the case if~$\PG$ is {\sl singular} (cf.~\cite{Voll/03b}).}. Thus the challenge to
compute the weight function~$w'([\Lambda'])$ is essentially reduced to
the problem of determining the $p$-adic valuation~$v_p(\det(M(\alpha^{d'})))$
(cf. equation~(\ref{weight-second-case})).

\section{Proofs}
\subsection{Proof of Theorem~\ref{theorem-flag-rational-function}} \label{subsection theorem rat function}
Choosing a common denominator for the sum~$F_n(q,{\bf X})$ we write 
$$F_n(q,{\bf X})=\frac{f_n(q,{\bf X})}{\prod_{i=1}^{n-1}(1-X_i)} $$
where 
\begin{eqnarray}
f_n(q,{\bf X})&=&\sum_{I\subseteq[{n-1}]}b_I(q)\prod_{i\in
  I}X_i\prod_{j\not\in I}(1-X_j)\label{numerator}
\\
&=&\sum_{I\subseteq[{n-1}]}c_I(q)\prod_{i\in I}X_i,\mbox{ say}.\nonumber
\end{eqnarray}
Then $$c_I(q)=\sum_{J\subseteq I\subseteq[{n-1}]}(-1)^{|I|-|J|}b_J(q)\in\mathbb{Z}[q].$$
For a subset~$I\subseteq[{n-1}]$ we define~$I^c:=[{n-1]}\setminus I$.
Theorem~\ref{theorem-flag-rational-function} will follow if we can prove
\begin{equation*}
c_I(q^{-1})=q^{-\binom{n}{2}}c_{I^c}(q)\quad\forall I\subseteq[{n-1}]
\end{equation*}
or, equivalently, if $c_I(q)=\sum_{k}a_{k,I}q^k$, that
\begin{equation}
a_{k,I}=a_{\binom{n}{2}-k,I^c}\quad\forall 0\leq k \leq \binom{n}{2},\;I\subseteq[{n-1}]\label{funeq-cIq-equivalent-formulation}.
\end{equation}
We shall prove~(\ref{funeq-cIq-equivalent-formulation}) by showing
that~$a_{k,I}$ enumerates the {\sl $k$-dimensional Schubert cell of
  type~$I$} in the complete flag variety~$\mcF_{[n-1]}$ and
that post-multiplication by the longest word induces a
$1-1$-correspondence between $k$-dimensional cells of type~$I$ and
$\left(\binom{n}{2}-k\right)$-dimensional cells of
type~${I^c}$. 

Firstly recall that the variety~$\mcF_{[{n-1}]}$ of complete flags in~$\Fq^n$ has a cell decomposition into {\sl Schubert cells}~$\Omega_w$, labelled by $w\in S_n$, the symmetric group on~$n$ letters (e.g.~\cite{Manivel/01}, Chapter~3). We write $w=(w_1\dots w_n)$ if $w(i)=w_i$. Let us represent a complete flag by an $n\times n$-matrix over~$\Fq$: its $r$-th member is generated by the first $r$ rows. For example~(\cite{Manivel/01}, p. 134) an element of~$\Omega_w$, $w=(365142)\in S_6$ has a unique matrix representative of the form
\begin{equation*}
\left(\begin{array}{cccccc}
*&*&1&0&0&0\\
*&*&0&*&*&1\\
*&*&0&*&1&0\\
1&0&0&0&0&0\\
0&*&0&1&0&0\\
0&1&0&0&0&0\\ 
\end{array}\right)
\end{equation*} 
Replacing the $*$'s by independent variables identifies the Schubert cell $\Omega_w$ with affine space over~$\Fq$ of dimension $l(w)$, where $l$ denotes the usual length function, the number of inversions of~$w$. The flag variety is now just the disjoint union of these affine spaces~$\Omega_w$. 

The expression~(\ref{numerator}) involves the cardinalities of the $2^{n-1}$ varieties~$\mathcal{F}_I$, $I\subset[{n-1}]$. The cell decomposition of~$\mcF_{[{n-1}]}$ will allow us to to accommodate all of the~$\mathcal{F}_I$ in one object by identifying them with certain unions of Schubert cells. 

To that end we define the {\sl type} $\nu(w)$ of a permutation~$w$ to
be the smallest subset~$I$ of $[{n-1}]$ such that the natural
{\sl surjection} $G/B\rightarrow G/P_I$ of complete flags onto flags
of type $I$ is a {\sl bijection} if restricted to
$\Omega_w$. Alternatively, given $w=(w_1\dots w_n)$, set
$$\nu(w):=\{i\in[{n-1}]|\; w_{i+1}<w_i\}.$$ So for~$w$ as in
our example above we have~$\nu(w)=\{2,3,5\}$, the type of the longest
word $w_0=(n\,n-1\dots2\,1)$ equals $[{n-1}]$, and the type of
the identity element is the empty set. We have a bijection of sets  
\begin{equation}
\coprod_{\nu(w)\subseteq I}\Omega_{w}\stackrel{1-1}{\longleftrightarrow} \mcF_I. \label{correspondence}
\end{equation}
From~(\ref{correspondence}) it follows immediately that
\begin{eqnarray}
\sum_{\nu(w)\subseteq I}|\Omega_{w}|&=&b_I(q)\quad \mbox{ and}\nonumber\\
\sum_{\nu(w)= I}|\Omega_{w}|&=&c_I(q).\nonumber
\end{eqnarray}
As the Schubert cells~$\Omega_{w}$ are identified with some affine
$l(w)$-space, their cardinalities are just powers of $q$ and
$a_{k,I}(q)$, the $k$-th coefficient of the polynomial~$c_I(q)$,
counts the number of Schubert cells of type~$I$ of dimension~$k$. Theorem~\ref{theorem-flag-rational-function} will follow from the following

\begin{proposition} \label{proposition_bijection}
Let $w_0=(n\,n-1\dots2\,1)\in S_n$ be the longest word. 
The bijection $w\mapsto w\,w_0$ induces bijections 
$$\{\Omega_w|\;\nu(w)=I\}\stackrel{1-1}{\longleftrightarrow}\{\Omega_w|\;\nu(w)={I^c}\}.$$
We have $l(w\,w_0)=\binom{n}{2}-l(w)$.%
\end{proposition}

\begin{proof} This follows easily from the definition of the type of a
  permutation and a comparison of the Schubert cells~$\Omega_{w}$ and
  $\Omega_{w\,w_0}$ (e.g. as sets of matrices).
\end{proof}

\subsection{Proof of Theorem~\ref{theorem-funeq} and
  Corollary~\ref{coro to main theorem}  }\label{subsection proof of theorem-funeq}
\subsubsection{The case $n_{\PG}(p)=0$}\label{subsubsection-n_0=0}

First we deal with the case that the Pfaffian hypersurface~$\PG$ has no $\Fp$-rational points. Thus
$\det(M(\alpha))$ is a $p$-adic unit for all
$\alpha\in\Zp^{d'}\setminus p\Zp^{d'}$ (i.e. for all column vectors of
matrices in~$Sl_{d'}(\Zp)$) and~(\ref{linear congruences}) is
equivalent to the single congruence
$$ {\bf \overline{g}} \equiv 0 \; \mod \; p^{\sum_{j=1}^l
  r_{i_j}}.$$ 
Hence 
\begin{equation}
|G_p:X(\Lambda')|=p^{d\sum_{j=1}^l
  r_{i_j}} \label{index abel}
\end{equation}
and
\begin{equation*}
w'([\Lambda'])=\sum_{i\in\nu([\Lambda'])}r_i(d+d'-i).\label{weight-first-case}
\end{equation*} 
Thus
\begin{eqnarray}
A(p,p^{-s})&=&\sum_{I\subseteq[{d'-1}]}A_I(p,p^{-s})\nonumber \\ 
&=&\sum_{I\subseteq[{d'-1}]}\sum_{{\bf r}_I>0}f(I,{\bf r}_I,p)\cdot p^{d\,{\sum r_{i}(d'-i)}-s\sum r_{i}(d+d'-i)} \nonumber\label{sum-no-points}\\
&\stackrel{*}{=}&\sum_{I\subseteq[d'-1]}\frac{b_I(p)}{p^{\dim \mcF_I}}\sum_{{\bf r}_I>0}p^{\sum r_{i}(d+i)(d'-i)-s\sum r_{i}(d+d'-i)}\nonumber \\
&\stackrel{**}{=}&\sum_{I\subseteq[d'-1]}{b_I(p^{-1})}\prod_{i\in I}\frac{X_i}{1-X_i}\nonumber \\
&=:&\sum_{I\subseteq[d'-1]}F_{d'}(I,p^{-1},{\bf X}_I)=F_{d'}(p^{-1},{\bf X})\nonumber
\end{eqnarray}
Here we used Lemma~\ref{lemma-lattices-of-given-type} for
equality~($*$) and equation~(\ref{funeq-varieties}) to
obtain~($**$). For~$1\leq i \leq d'-1$, we made the substitutions
\begin{equation}
X_i:=p^{(d+i)(d'-i)-s(d+d'-i)}.\label{subsitutions}
\end{equation} 
Here ${\bf X}_I$ stands for $(X_i)_{i\in I}$, ${\bf r}_I>0$ for~${\bf
  r}_I\in\mathbb{N}_{>0}^{l}$ and~$\sum$ for~$\sum_{i\in
  I}$. Theorem~\ref{theorem-funeq} and its Corollary~\ref{coro to main
  theorem} follow now from
  Theorem~\ref{theorem-flag-rational-function} together with
  Corollary~\ref{coro to gss-lemma}.

\subsubsection{The case $n_{\PG}(p)>0$} 
In this case $F_{d'}(p^{-1},{\bf X})$ fails to represent the
generating function~$A(p,p^{-s})$, as~(\ref{index abel}) will not hold in
general. We shall see, however, that the two rational functions agree
`almost everywhere' and we will show how to decompose them into summands in
a geometrically meaningful way to see  exactly where and how they differ.
As we assume that the Pfaffian contains no lines we have
(by~(\ref{equiv congruence}))
\begin{equation}
A_I(p,p^{-s})=F_{d'}(I,p^{-1},{\bf X}_I)\mbox{ if }d'-1\,\not\in\,I.\label{agreement-1}
\end{equation}
For a point~${x}\in{\mathbb{P}^{d'-1}(\Fp)}$ and ${\bf X}':=(X_i)_{i\in[{d'-2}]}$ we set
\begin{eqnarray}
A(x,p,p^{-s})&:=&\sum_{\substack{d'-1\in\nu([\Lambda']) \\P([\Lambda'])=x}}p^{w([\Lambda'])\cdot d-sw'([\Lambda'])}\label{A_x-expression}\\
F_{d',0}(p^{-1},{\bf X}')&:=&F_{d'-1}(p^{-1},{\bf X}')p^{-(d'-1)}\frac{X_{d'-1}}{1-X_{d'-1}}\label{F_d'-expression}
\end{eqnarray}
The rational function~(\ref{A_x-expression}) might be thought of as the
generating function~(\ref{A}) with summation restricted to the maximal lattices
lifting a fixed point~$x\in\mathbb{P}^{d'-1}(\Fp)$. It agrees
with~(\ref{F_d'-expression}) if and only if $\det(M(x))\not=0\in\Fp$,
i.e. if~$x\not\in \overline{\PG}$:
\begin{equation}
A(x,p,p^{-s})=F_{d',0}(p^{-1},{\bf X}')\mbox{ if }x\not\in \overline{\PG}.\label{agreement-2}
\end{equation}
If we write
\begin{eqnarray*}
A(p,p^{-s})&=&\sum_{d'-1\not\in I}A_I(p,p^{-s})+\sum_{x\in\mathbb{P}^{d'-1}(\Fp)}A(x,p,p^{-s})\\
F_{d'}(p^{-1},{\bf X})&=&\sum_{d'-1\not\in I}F_{d'}(I,p^{-1},{\bf X}_I)+\binom{d'}{1}_pF_{d',0}(p^{-1},{\bf X}').
\end{eqnarray*}
we see that in order to prove Theorem~\ref{theorem-funeq} we are left
with the challenge to prove
\begin{proposition}\label{proposition-differences} Let
  $x\in\overline{\PG}$. Then
\begin{eqnarray}
A_{d'}(x,p,p^{-s})-F_{d',0}(p^{-1},{\bf
  X}')&=&F_{d'-1}(p^{-1},{\bf X}')p^{-(d'-1)}\frac{pY-X_{d'-1}}{(1-X_{d'-1})(1-Y)},\label{function-differences}
\end{eqnarray}
where
\begin{eqnarray*}
X_{d'-1}&:=&p^{d+d'-1}T^{d+1}\quad(\mbox{in accordance with }(\ref{subsitutions}))\mbox{ and }\\
Y&:=&p^{d+d'-2}T^{d-1}.
\end{eqnarray*}
In particular, the function $A_{d'}(x,p,p^{-s})-F_{d',0}(p^{-1},{\bf
  X}')$ is independent of~$x$.
\end{proposition}
Indeed, together with equations~(\ref{agreement-1}) and~(\ref{agreement-2})
Proposition~\ref{proposition-differences} clearly implies
\begin{eqnarray}
A(p,p^{-s})&=&F_{d'}(p^{-1},{\bf
  X})+\left(A(p,p^{-s})-F_{d'}(p^{-1},{\bf X})\right)\nonumber \\
&=&F_{d'}(p^{-1},{\bf
  X})+n_{\PG}(p)\left(A_{d'}(x,p,p^{-s})-F_{d',0}(p^{-1},{\bf
  X}')\right)\nonumber\\
&=&F_{d'}(p^{-1},{\bf
  X})+n_{\PG}(p)F_{d'-1}(p^{-1},{\bf X}')p^{-(d'-1)}\frac{pY-X_{d'-1}}{(1-X_{d'-1})(1-Y)}.\label{A(p,p^{-s}-expression}
\end{eqnarray}

\begin{proofnodot} (of Proposition~\ref{proposition-differences}) Let~$J_2$ denote the
  matrix~$\left(\begin{array}{cc}0&1\\-1&0\end{array}\right)$. Locally
  around any of the $n_{\PG}(p)$
points of the Pfaffian mod~$p$, the admissibility
  conditions~(\ref{linear congruences}) look like
\begin{eqnarray}
\bfg\cdot\mbox{diag}\left(\left(\begin{array}{cc}0&x_1\\-x_1&0\end{array}\right),J_2,\dots,J_2\right)&\equiv& 0 \mod p^{\sum_{j=1}^{l} r_{i_j}}  \label{admiss-condition-1}\\
\bfg&\equiv&0 \mod p^{\sum_{j=1}^{l-1}r_{i_j}},\nonumber
\end{eqnarray}
where ${\bf x }=(x_1:\dots:x_{d'})\in\mathbb{P}^{d'-1}(\Zp/(p^{r_{d'-1}}))$, ${\bf x}\equiv (0:1:\dots:1)\mod p$ 
and we have
\begin{equation}
w'([\Lambda'])=\sum_{i\in\nu([\Lambda'])}r_{i}(d+d'-i)-2\min\{r_{d'-1},v_p(x_1)\}.\label{weight-second-case}
\end{equation}
Therefore
\begin{eqnarray}
A_{d'}(x,p,p^{-s})&=&B_0(p,p^{-s})F_{d'-1}(p^{-1},{\bf X}')\label{A_d'-expression}
\end{eqnarray} say, where
\begin{equation}
B_0(p,p^{-s}):=\sum_{\substack{r_{d'-1}>0\\({\bf x
    })\in\mathbb{P}^{d'-1}(\Zp/(p^{r_{d'-1}}))\\{\bf x}\equiv
    (0:1:\dots:1)\mod
    p}} p^{dr_{d'-1}-s(1+d-2\min\{r_{d'-1},v_p(x_1)\})} \label{B_0}
\end{equation}
We see that in the present case the weight function~$w'([\Lambda'])$ given in (\ref{weight-second-case})  depends on more than just the lattice's
type. To get an explicit expression for the function~$B_0(p,p^{-s})$, we
must find a way to eliminate the term `min' in the
sum~(\ref{B_0}). In other words, we must answer the following
question: Given~$(a,b)\in\mathbb{N}^2_{>0}$, how many of the
lifts~$(x_1:\dots:x_{d'})$ of the point $(0:1:\dots:1)\in\mathbb{P}^{d'-1}(\Fp)$ to points mod~$p^a$ are there with~$v_p(x_1)=b$?
\begin{lemma} For $(a,b)\in N:=\{(x,y)\in\mathbb{N}^2_{>0}|\,x\geq y \geq 1\}$ let
$$\lambda(a,b):=\{{\bf x}\in\mathbb{P}^{d'-1}(\Zp/(p^a))|\,{\bf x}\equiv (0:1:\dots:1)\mod p,\,v_p(x_1)=b\}. $$
Then
\begin{eqnarray*}
\lambda(a,b)&=&\left\{\begin{array}{ll}p^{(d'-2)(a-1)} & \mbox{ if }(a,b)\in\Delta,\\
                                       p^{(d'-2)(a-1)+a-b}(1-p^{-1})&\mbox{ if }(a,b)\in N\setminus\Delta, \end{array}\right. 
\end{eqnarray*}
where $\Delta:=\{(x,y)\in N|\,x=y\}$.
\end{lemma}
\begin{proof}
This is easy to check in an affine chart.
\end{proof}

\noindent Thus
\begin{eqnarray}
B_0(p,p^{-s})&=&\sum_{(a,b)\in\Delta}p^{(d'-2)(a-1)+ad-s(d-1)a}+\nonumber\\&&\quad
(1-p^{-1})\sum_{(a,b)\in N\setminus\Delta}p^{(d'-2)(a-1)+a-b+ad-s((d+1)a-2b)}\nonumber\\
&=&p^{-(d'-2)}\left(\frac{Y}{1-Y}\left(1+(1-p^{-1})\frac{X_{d'-1}}{1-{X_{d'-1}}}\right)\right)\nonumber\\
&=&p^{-(d'-1)}\frac{Y(p-X_{d'-1})}{(1-Y)(1-X_{d'-1})},\label{B_0-expression}
\end{eqnarray}
where $X_{d'-1}$ and $Y$ are defined as in the statement of
Proposition~\ref{proposition-differences}, which now follows from routine computations combining the identities~(\ref{F_d'-expression}),~(\ref{A_d'-expression}) and~(\ref{B_0-expression}).
\end{proofnodot}

With equation~(\ref{A(p,p^{-s}-expression}) we have given an explicit
formula for the generating function~$A(p,p^{-s})$, which, by
Lemma~\ref{gss-lemma}, is tantamount to the local normal zeta function,
completing the proof of Theorem~\ref{theorem-funeq}. The functional equation also follows swiftly
from~(\ref{A(p,p^{-s}-expression}), Theorem~\ref{theorem-flag-rational-function} and Corollary~\ref{coro to gss-lemma} if the Pfaffian hypersurface~$\PG$ is absolutely
irreducible. Indeed, let $V$ be {\sl any} non-singular, absolutely irreducible
projective variety over~$\Fp$ of dimension~$n$. If~$b_{V,e}$,
$e\geq1$, denotes the number of $\mathbb{F}_{p^e}$-rational points of
$V$ it is a well-known consequence of the {\sl rationality} of the Weil zeta function
$$Z_V(u)=\exp\left(\sum_{e=1}^\infty \frac{b_{V,e}u^e}{e}\right)$$
that there are complex numbers $\beta_{r,j}$, $r=0,\dots,2n$, $j=1,\dots,B_r$, $B_r\in \N$, such that 
$$b_{V,e}=\sum_{r=0}^{2n}(-1)^r\sum_{j=1}^{t_r}\beta^e_{r,j},$$
and that the function
\begin{eqnarray*}
\N_{\geq0}&\rightarrow&\N \\
e&\mapsto&b_{V,e}
\end{eqnarray*}
 has a {\sl unique} extension to $\Z$ (cf~\cite{DenefMeuser/91}, Lemma~2). 
The {\sl functional equation} of the Weil zeta function
$$Z_V(1/p^nu)=\pm (p^{n/2}u)^\chi Z_V(u),$$
where $\chi = \sum_{i=1}^{2n} (-1)^iB_i$, implies the $1-1$-correspondences
$$\left\{\frac{p^n}{\beta_{r,j}}|\;1\leq j \leq
B_r\right\}\stackrel{1-1}{\longleftrightarrow}\left\{\beta_{2n-r,i}|\;1
\leq i \leq B_{2n-r}\right\}$$
for $0\leq j \leq 2n$ (cf~\cite{Igusa/00}, p.~213). This gives 
$$b_{V,-e}=p^{-en}\,b_{V,e}$$ formally\footnote{An instance of which
  we have also seen in~(\ref{funeq-varieties}) where $b_{V,e}$ was
  given as a polynomial in~$p^e$.}. Corollary~\ref{coro to main
  theorem} now follows immediately if we set $V=\PG$, $e=1$, $n=d'-2$,
$n_{\PG}(p^e)=b_{V,e}$.

\subsection{Proof of Theorem~\ref{theorem-grenham}}\label{subsection grenham}
The proof of Theorem~\ref{theorem-grenham} would have been presented
in Section~\ref{subsubsection-n_0=0} had this not interrupted the proof
of Theorem~\ref{theorem-funeq}. Here~$d=n$ and~$d'=n-1$. The essential
observation is that again the weight function~$w'([\Lambda'])$ only
depends on the lattice's {\sl type}. It was indeed explicitly
calculated in~\cite{Voll/02}, Chapter~5.2, as
$$w'([\Lambda'])=\sum_{i\in\nu([\Lambda'])}r_i(2(d'-i)+1).$$ This allows us to write

\begin{eqnarray*}
A(p,p^{-s})&=&\sum_{I\subseteq[{d'-1}]}A_I(p,p^{-s})\nonumber \\ 
&=&\sum_{I\subseteq[{d'-1}]}\sum_{{\bf r}_I>0}f(I,{\bf r}_I,p)\cdot p^{d\,{\sum r_{i}(d'-i)}}T^{\sum r_{i}(2(d'-i)+1)} \label{grenham}\\
&=&\sum_{I\subseteq[{d'-1}]}\frac{b_I(p)}{p^{\dim \mcF_I}}\sum_{{\bf r}_I>0}p^{\sum r_{i}(d+i)(d'-i)}T^{\sum r_{i}(2(d'-i)+1)}\nonumber \\
&=&F_{d'}(p^{-1},\tilde{{\bf X}}).\nonumber
\end{eqnarray*}
Here, for~$1\leq i \leq d'-1$, we made the substitutions
\begin{equation*}
\tilde{X}_i:=p^{(d+i)(d'-i)}T^{2(d'-i)+1}.
\end{equation*} 
Again ${\bf r}_I>0$ stands for~${\bf r}_I\in\mathbb{N}_{>0}^{l}$
and~$\sum$ for~$\sum_{i\in I}$. The result follows
from Theorem~\ref{theorem-flag-rational-function} and Corollary~\ref{coro to gss-lemma}.

\newpage
%
%


\def\cqfd{\kern 2truemm\unskip\penalty 500\vrule height
4pt depth 0pt width 4pt\medbreak} \def\carre{\vrule height
5pt depth 0pt width 5pt} \def\virg{\raise .4ex\hbox{,}}
\def\pf{\mathop{\rm Pf}\nolimits}
\def\og{\leavevmode\raise.3ex\hbox{$\scriptscriptstyle\langle\!
\langle$}}
\def\fg{\leavevmode\raise.3ex\hbox{$\scriptscriptstyle\rangle\!
\rangle$}}
\input amssym.def
\def\P{\mathbb{P}}
\font\san=cmssdc10
\def\ext{\hbox{\san \char3}}
\frenchspacing
\appendix
\section*{}

\centerline{\sl Appendix}
\bigskip
\centerline{\bf Lines on pfaffian hypersurfaces}
\smallskip
\centerline{Arnaud BEAUVILLE}\bigskip
The aim of this Appendix is to prove
 that a general pfaffian hypersurface of degree $r>2n-3$ in $\P^n$ contains no lines (Proposition~\ref{beauville-proposition}). By a
simple dimension count (see Corollary 2 below), it suffices to
show that the variety of lines contained in the universal pfaffian
hypersurface (that is, the hypersurface of degenerate forms in the
space of all skew-symmetric forms
 on a given vector space)  has the
expected dimension. We will deduce this from an explicit
description of the  pencils of degenerate skew-symmetric forms,
which is the content of the Proposition below.\medskip
\par We work over an algebraically closed field $k$. We
will need an elementary lemma:

\begin{lemma}
Given a
pencil of  skew-symmetric forms on a $n$-dimensional
vector space, there exists a subspace of dimension $[\frac{n+1}{2}]$
 which is isotropic for all forms of the pencil.
\end{lemma}
\begin{proof}
By induction on $n$, the cases $n=0$ and
$n=1$ being trivial. Let
$\varphi +t\psi$ be our pencil; we can assume that  $\varphi $ is
degenerate. Let $D$ be a line contained in the kernel of $\varphi
$, and let $D^{\perp}$ be its orthogonal with respect to~$\psi$.
Then
$\varphi $ and $\psi$ induce skew-symmetric forms $\bar\varphi
$ and $\bar\psi$ on
$D^{\perp}/D$; by the induction hypothesis there exists a
subspace   of dimension $[\frac{n-1}{2}]$ in $D^{\perp}/D$ which is isotropic for
$\bar\varphi$ and $\bar\psi$. The pull-back of this subspace in
$D^{\perp}$ has dimension $[\frac{n+1}{2}]$ and is  isotropic for
$\varphi$ and $\psi$.
\end{proof}

The following result must be well-known, but I have not been
able to find a  reference:

\begin{proposition}\label{beauville proposition}
Let
$V$ be a vector space of dimension
$2r$, and  $(\varphi _t)_{t\in \P^1}$  a pencil of {\rm degenerate}
skew-symmetric forms on $V$. There exists a subspace
$L\in V$ of dimension $r+1$ which is isotropic for
$\varphi _t$ for all $t\in \P^1$.
\end{proposition}

\begin{proof}
Again we prove the Proposition by induction on
$r$, the case $r=1$ being trivial. The associated maps
$\Phi_t:V\rightarrow V^*$ form a pencil of singular
linear maps. By a classical result in linear algebra (see~\cite{Gantmacher/84}, chap.
XII, thm. 4), there exist subspaces $K\in V$ and $L'\in V^*$, with
$\dim K=\dim L'+1$, such that  $\Phi_t(K)\in L'$ for all $t$;
equivalently, there exist subspaces $K$ and $L$ of $V$, with
$\dim K+\dim L=2r+1$, which are orthogonal for each
$\varphi_t$. Replacing $(K,L)$ by $(K\cap L,K+L)$ we may
assume $K\in L$; the pencil $(\varphi_t)$ restricted to $L$ is
 singular on $K$, hence induces a pencil  $(\bar\varphi_t)$
on $L/K$. Put $\dim K=p$, so that  $\dim (L/K)=2r+1-2p$.
By the above lemma  there is a subspace of $L/K$, of
dimension
$r+1-p$, which is isotropic for each $\bar\varphi_t$. Its pull
back in $L$ has dimension $r+1$ and is isotropic for each
$\varphi_t$.
\end{proof}

Let us give a few consequences of Proposition~\ref{beauville proposition}. We
keep our vector space~$V$ of dimension~$2r$; we denote by ${\cal S}_r$  the space of skew-symmetric
forms on $V$, and  by
$\Xr$ the hypersurface of
degenerate forms in $\P({\cal S}_r)$.

\begin{corollary} \label{beauville corollary 1}
The variety of lines contained
in $\Xr$ is irreducible, of codimension $r+1$ in the
Grassmannian of lines of $\P({\cal S}_r)$.
\end{corollary}

\begin{proof} The $(r+1)$-planes  of $V$
are parametrized by a Grassmannian ${\cal G}$ of dimension
$r^2-1$. For  such a plane $L$ the space ${\cal S}_{r,L}$ of forms
$\varphi\in{\cal S}_r$ vanishing on $L$ has dimension
$$\dim{\cal S}_{r,L}=\dim \ext^2V^* -\dim
\ext^2L^*=r(2r-1)-\frac{r(r+1)}{2}=\frac{3r(r-1)}{2}\ .$$
Let ${\cal P}$ be the Grassmannian of lines in $\P({\cal S}_r)$
(that is, the variety of pencils of skew-symmetric forms). Consider
the locus $Z\in {\cal P}\times {\cal G}$ of pairs $(\ell ,L)$ with $\ell \in {\cal S}_{r,L}$.  The
projection $Z\rightarrow {\cal G}$ is a smooth fibration; its
fibre  above a point $L\in {\cal G}$ is the
Grassmannian
 of lines in $\P({\cal S}_{r,L})$, which has dimension
$2\dim{\cal S}_{r,L}-4$. Thus~$Z$ is smooth, irreducible, of
dimension $r^2-1+2\dim{\cal
S}_L-4=4r^2-3r-5.$

Let ${\cal P}_{\rm sing}$ be  the subvariety of
${\cal P}$ consisting of lines contained in $\Xr$ (that
is, the subvariety of singular pencils). The content of
Proposition~\ref{beauville proposition} is that ${\cal P}_{\rm sing}$ is the image of
$Z$ under the projection to ${\cal P}$. Thus ${\cal P}_{\rm sing}$ is irreducible, of
dimension $\le 4r^2-3r-5$, or equivalently,  since
$\dim{\cal P}= 2\dim{\cal S}_r-4=4r^2-2r-4 $, of codimension $\ge
r+1$. On the other hand, ${\cal P}_{\rm sing}$ is defined locally
by $(r+1)$ equations in ${\cal P}$, given by the coefficients of
the polynomial $\pf(\varphi _t)$ of degree~$r$. The Corollary
follows.
\end{proof}

Observe that  $r+1$ is  the number of conditions that the requirement to
contain a given line imposes on a hypersurface of degree $r$
in projective space. In other words, Corollary~\ref{beauville corollary 1} says  that
the hypersurface $\Xr$ behaves like a general
hypersurface of degree $r$ as far as the dimension of its variety
of lines is concerned.
\medskip
\par Let $L$ be a   vector space, of dimension $n+1$, and
$\ell =(\ell _{ij})$ a
$(2r\times 2r)$-skew-symmetric matrix of linear forms on $L$.
The hypersurface $X_\ell $ in $\P(L)\ (=\P^n)$ defined by
$\pf(\ell _{ij})=0$ is called a {\it pfaffian
hypersurface}. It is defined by the equation $\pf(\ell _{ij})=0$, of
degree~$r$.

\begin{corollary}
If $r>2n-3$ and  the forms
$\ell _{ij}$ are general enough, $X_\ell $ contains no lines.
\end{corollary}
\begin{proof}
The matrix $(\ell _{ij})$
defines a linear map $u:L\rightarrow {\cal S}_r$, which is injective
 when the forms $\ell _{ij}$ are general enough (observe that
$\dim L<\dim {\cal S}_r $). Thus we can identify $L$ to its image in
${\cal S}_r$, and $X_\ell $ to the hypersurface $\Xr\cap
\P(L)$ in~$\P(L)$.
\par  Let $G$ be the Grassmann
variety of
 $(n+1)$-dimensional vector subspaces of~${\cal
S}$, and
$F$ the variety of lines contained in~$\Xr$. Consider
 the incidence variety $Z\in
F\times G$ of pairs $(\ell ,L)$ with $\ell \in \P(L)$. The fibre of the
projection $Z\rightarrow G$ at a point $L\in G$ is the variety of
lines contained in $\Xr\cap\P(L)=X_\ell $.
\par Put $N:=\dim {\cal S}_r$. We have
$\dim F=2N-4-(r+1)$ by Corollary 1; the projection
$Z\rightarrow F$ is a fibration  of relative dimension
$(n-1)(N-n-1)$. This gives $\dim Z=2N-4-(r+1)+(n-1)(N-n-1)$,
while  $\dim G=(n+1)(N-n-1)$. Thus
$$\dim Z - \dim G= 2n-3-r<0\ ,$$hence the general
fibre of the projection
$Z\rightarrow G$ is empty.
\end{proof}

Note that `$(\ell _{ij})$ general enough' means
`for $(\ell _{ij})$ in a certain Zariski open subset of
$(L^*)^N$'.
In particular, suppose
that our vector space $L$ comes
from a vector space $L_0$  over  an infinite subfield $k_0$
of $k$; then {\it
 the matrices $(\ell _{ij})\in (L_0^*)^N$
 such that $X_\ell $ contains no lines are
Zariski dense in the parameter space $(L^*)^N$ for} $r>2n-3$.







\bibliographystyle{amsplain} 
\bibliography{thebibliography}
\vskip1cm
\begin{tabular}{lll}
Christopher Voll &\hfill& Arnaud Beauville \\
Mathematical Institute  &\hfill& Institut Universitaire de France \&\\
24 - 29 St. Giles' &\hfill& Laboratoire J.-A. Dieudonn\'{e}\\
Oxford  OX1 3LB &\hfill& UMR 6621 du CNRS\\
United Kingdom &\hfill& Universit\'{e} de Nice \\
&\hfill& Parc Valrose \\
&\hfill& 06108 Nice Cedex 2 \\
&\hfill& France\\
\end{tabular}

\end{document}